\documentclass[11pt]{article}

\def\mathbb#1{{\bf #1}}

\usepackage{a4wide,epsfig,times,amsmath,amstext,amssymb}

\newcommand{\qed}{\square} 
\newtheorem{theorem}{Theorem}[section]
\newtheorem{definition}[theorem]{Definition}

\newtheorem{prop}[theorem]{Proposition}
\newtheorem{lemma}[theorem]{Lemma}
\newtheorem{corollary}[theorem]{Corollary}
\newtheorem{example}[theorem]{Example}

\newenvironment{proof}{\noindent {\bf Proof} }{$\qed$}

\newcommand{\height}{\operatorname{height}}

\newcommand\ra{\rightarrow}

\newcommand\proper[1]{{#1}^\circ}
\newcommand\atom[1]{\operatorname{atom}({#1})}
\newcommand\ceiling[1]{\lceil{#1}\rceil}
\newcommand\floor[1]{\lfloor{#1}\rfloor}
\newcommand\supp{\operatorname{supp}}
\newcommand\lcm{\operatorname{lcm}}
\newcommand\link{\operatorname{link}}

\newcommand\mi[1]{\operatorname{mi}({#1})}
\newcommand\reg{\operatorname{reg}}
\newcommand\codim{\operatorname{codim}}
\newcommand\rank{\operatorname{rank}}
\newcommand\length{\operatorname{length}}
\newcommand\pol{\operatorname{pol}}
\newcommand\scarf{\operatorname{Sc}}

\newcommand\width{\operatorname{width}}
\newcommand\pd{\operatorname{pd}}

\begin{document}

\begin{center} 
{\bf \Large{Minimal monomial ideals and linear resolutions}} \\
Jeffry Phan \\
Columbia University
\end{center}

\begin{abstract}
A minimal monomial ideal is the combinatorially simplest monomial ideal whose lcm-lattice equals a given finite atomic lattice $\hat{L}$.  The minimal ideal inherits  many nice properties of any ideal $I$ whose lcm-lattice also equals $\hat{L}$, e.g. Cohen-Macaulayness and the dual property of having a linear resolution.  Conversely, any ideal having a linear resolution is shown to be, essentially, minimal.
\end{abstract}

\setcounter{section}{0}

\begin{section}{Introduction}\label{introduction section}
We introduce and study minimal monomial ideals; these are the simplest ideals whose lcm-lattice equals a given finite atomic lattice $\hat{L}$.  Recall that if $I\subset k[y_1,\ldots,y_n]$ is a monomial ideal minimally generated by $f_1,\ldots, f_r$, the lcm-lattice $LCM(I)$ of $I$ is the set of all least common multiples of the $f_i$, partially ordered by divisibility.  $LCM(I)$ is a finite atomic lattice; Gasharov, Peeva, and Welker proved $LCM(I)$ computes the Betti numbers of $I$ (and more, see \cite{GPW}).  They also showed that if $LCM(I)\cong LCM(J)$, then $I$ and $J$ have equivalent resolutions up to a relabeling process they make precise (see also the relabeling used in the deformation of exponents of~\cite{BPS} and~\cite{MSY}).  Our construction shows every finite atomic lattice is the lcm-lattice of some monomial ideal.  

Although the minimal free resolution of a monomial ideal is easily computed, it is unknown how to theoretically describe the maps between the free modules.  When $I$ has a linear resolution, Reiner and Welker constructed these maps in~\cite{Reiner_Welker}.  Bayer, Peeva, and Sturmfels showed the maps can be thought of as the boundary maps of a simplicial complex when $I$ is a generic monomial ideal (\cite{BPS}) and this idea has been studied in several papers (\cite{BS}, \cite{MSY}, \cite{Batzies_Welker}).  It is natural to ask if $LCM(I)$ can be used to describe the maps in the minimal resolution.   This question is still widely open, but minimal monomial ideals are one possible avenue of attack: if one can construct the maps for minimal ideals, then up to relabeling the problem is solved for all monomial ideals.

In Section~\ref{constructions section} we construct minimal squarefree ideals and prove they are universal in the sense of Theorem~\ref{universal property theorem}.  One corollary is that every simplicial complex $\Delta$ is the Scarf complex of some squarefree monomial ideal (\cite{BPS},\cite{BS}).  If $\Delta$ is acyclic, then $\Delta$ supports the minimal resolution of that ideal.

The subset $\mi{L}\subset{L}$ of meet-irreducible elements of $L$ plays an essential role throughout this paper.  In Section~\ref{distributive section}, we study how $\hat{L}$ is related to the distributive lattice generated by $\mi{L}$.    Closely related to minimal ideals are the the Hibi ring of $\mi{L}$ and the ideals ${\cal{H}_{\hat{L}}}$ studied in the series of papers by Herzog, Hibi, and Zheng (\cite{HH1},\cite{HH2},\cite{HHZ},\cite{HHZ1}); exploration of the connections are left for a future paper.  

In Section~\ref{depolarization section}, we study nonsquarefree minimal ideals.  As a corollary, we are able to prove the width of $\mi{L}$ bounds the projective dimension of $I$.  

In Section~\ref{inherited properties section}, we show that many nice properties of $I$ pass to the minimal ideal constructed from $LCM(I)$.  In particular, if $I$ has a linear resolution, then so does the minimal ideal.  We characterize those $\hat{L}=LCM(I)$ for which this is possible; a very similar result (by a theorem of Eagon and Reiner) was proved by Yuzvinsky in the setting of rings of sheaves on posets (\cite{Yuzvinsky}).  Finally, we show that ideals with linear resolutions are minimal up to a common divisor of the generators.
\end{section}

\begin{section}{Squarefree Constructions}\label{constructions section}
Let $\hat{L}$ be a finite atomic lattice with proper part $L = \hat{L} - \{\hat{0},\hat{1}\}$ and let $\mi{L}$ be the subposet of meet-irreducible elements of $L$.  Fix a field $k$.  We will construct the simplest possible squarefree monomial ideal $M(L)\subset k[L]:=k[x_l :\, l\in \mi{L}]$ whose lcm-lattice equals $\hat{L}$.  

For $a\in \hat{L}$ denote by $\floor{a}$ and $\ceiling{a}$ the order ideal $\floor{a} := \{ b\in \hat{L}:\, b\leq a \}$ and the order filter $\ceiling{a}:=\{ b\in \hat{L} :\, a \leq b\}$ generated by $a$.  Define $x(a):= \prod_{l} x_l \in k[L]$ where the product is over all $l\in \mi{L} - \ceiling{a}$.  Partially order monomials by divisibility.

Every $a\in \hat{L}$ equals the meet of those $l\in \mi{L}$ satisfying $a\leq l$.  This implies $a\leq b$ in $\hat{L}$ if and only if $\mi{L} \cap \ceiling{b} \subset \mi{L} \cap \ceiling{a}$ if and only if $\mi{L} - \ceiling{a} \subset \mi{L} - \ceiling{b}$.  Therefore
\begin{equation}\label{separation equation}
a\leq b \qquad \text{if and only if} \qquad x(a)\leq x(b).
\end{equation}

\begin{theorem}\label{main construction theorem}
Let $a_1,\ldots, a_r$ be the atoms of $\hat{L}$.  Define the squarefree monomial ideal $M(L)\subset k[L]$ by
$$M(L):=(x(a_1),\ldots,x(a_r)).$$
Then $LCM(M(L))\cong \hat{L}$.
\end{theorem}
\begin{proof}
Set $L'=LCM(M)$.  Let $b\in \hat{L}$ and let $\supp(b):=\{i:\, a_i \leq b\}$ be the support of $b$.  Then $\mi{L}-\ceiling{b}$ equals the union of the $\mi{L}-\ceiling{a_i}$, $i \in \supp(b)$, and therefore $x(b) = \lcm(x(a_i):\, i \in \supp(b))$.  It follows that $x(b\vee c) = \lcm(x(b),x(c))$ for all $b,c\in \hat{L}$.  This means the map $x:\hat{L}\ra L'$, $a\mapsto x(a)$, preserves joins and is surjective.  $x$ is determined by supports and $\supp(b\wedge c) =\supp(b)\cap \supp(c)$ so $x$ also preserves meets.  $x$ is injective because $b\nleq c$ implies $x(b)\nleq x(c)$, by Equation~\ref{separation equation}.  Therefore $x$ is a lattice isomorphism.
\end{proof}

\begin{example}\rm
The two examples in Figure 1 illustrate Theorem~\ref{main construction theorem}.  Each $k\in \mi{L}$ has been colored light grey and labeled with the variable $k$ instead of $x_k$.  The hats $\, \hat{\phantom{}} \,$ are meant to suggest that every element in $\floor{k}$ is missing the variable $k$.  The lefthand figure is the the lcm-lattice of the ideal $(bd,cd,ac)\subset k[a,b,c,d]$.  The righthand figure is the lcm-lattice of $(befg,dfg,ceg,acd,bdef)\subset k[a, \ldots, f]$; in this example the labels of non-atoms have been omitted for clarity.
\begin{figure}[h]
\centering
\mbox{\epsfig{file=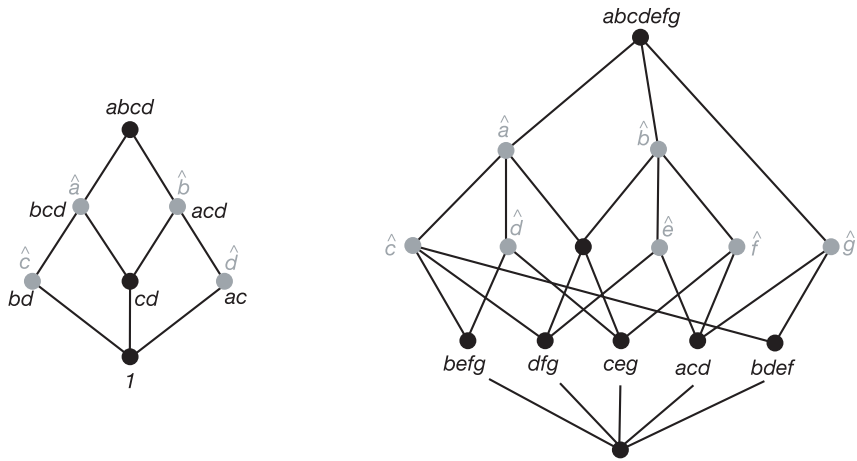, height = 6cm}}
\caption{Two atomic lattices and the corresponding minimal squarefree ideals.}
\end{figure}
\end{example}

\begin{example}\rm\label{cox example}
Let $\partial P$ be the boundary complex of the convex polytope $P$ and let $L$ be set of faces of $\partial P$ ordered by inclusion.  The atoms $a$ of $\hat{L}$ are the vertices of $\partial P$ and the coatoms of $\hat{L}$ are the facets of $\partial P$.  Because $\hat{L}$ is coatomic, $x(a)$ is the product of those variables corresponding to facets not containing the vertex $a$.  Hence $M(L)$ is the Cox irrelevant ideal of the toric variety defined by $P$ (\cite{Cox}).
\end{example}

\begin{example}\rm
When $L$ is a geometric lattice, $M(L)$ coincides with the ideal constructed by Irena Peeva in~\cite{Peeva}.  In that paper, Peeva proved the minimal free resolution of $I$ can be explicitly described using either broken-circuit complexes or Orlik-Solomon complexes.  Isabella Novik proved the same can be done using rooted complexes (\cite{Novik}).
\end{example}

Recall the Taylor complex $T$ of $M(L) = (x(a_1), \ldots , x(a_r))$ is the simplex on $[r]:= \{1,\ldots,r\}$ with each face $F\in T$ labeled by the degree $x(F):= \lcm \{ x(a_i):\, i \in F\}$.  The Scarf complex $\scarf{M(L)}\subset T$ is the subcomplex of faces that have unique degree.  Assume $\scarf{M(L)}$ is not a simplex.  The face poset of $\scarf{M(L)}$ is then naturally a subposet of $\hat{L}$ and therefore ideals with isomorphic lcm-lattices must have isomorphic Scarf complexes.  $\scarf{M(L)}$ supports a minimal resolution of $M(L)$ if and only if each Betti degree of $M(L)$ is the degree of some face of $\scarf{M(L)}$ if and only if for all $b\in \hat{L}$ either the closed interval $[\hat{0},b]_{\hat{L}}$ is boolean or the open interval $(\hat{0},b)_{\hat{L}}$ is acyclic.  The last equivalence follows from the next theorem, which appears as Theorem 2.1 in~\cite{GPW}; we call a multidegree $b$ for which $\beta_i(k[L]/M(L),b)\neq 0$ a {\em Betti degree} of $k[L]/M(L)$.  Note that we identify the interval $(\hat{0},b)_{\hat{L}}$ topologically with its order complex.

\begin{theorem}{(\cite{GPW})}\label{betti degree theorem}
For $i\geq 1$ and $b\in \hat{L}$ we have $\beta_i(k[L]/M(L),b) = \dim_k \tilde{H}_{i-2} \left( (\hat{0}, b)_{\hat{L}}; k \right)$, where $\beta_i(k[L]/M(L),b)$ is the $i$th Betti number in multidegree $b$.
\end{theorem}

\begin{prop}\label{scarf prop}
Every simplicial complex $\Delta$ not equal to the boundary of a simplex is the Scarf complex of some squarefree monomial ideal.  If $\Delta$ is acyclic, then $\Delta$ supports the minimal resolution of that ideal. 
\end{prop}
\begin{proof}
If $\Delta$ is a simplex, then $\Delta$ is the Scarf complex of the irrelevant ideal and supports the minimal resolution of that ideal.  Let $\hat{L}$ be the face poset of $\Delta$.  We can assume that $\Delta$ is not a simplex or, equivalently, that $\hat{L}$ is atomic.  Let $M(L)$ be as in Theorem~\ref{main construction theorem}, so that $LCM(M(L))=\hat{L}$.  It's easy to see $\Delta$ is the Scarf complex of $M(L)$.  $\hat{1}\in \hat{L}$ is the least common multiple of all the generators.  $\hat{1}$ is not a betti degree if and only if $(\hat{0},\hat{1})_{\hat{L}}$ is acyclic if and only if $\Delta$ is acyclic; in this case, the betti degrees are concentrated in $\Delta$ so $\Delta$ supports the minimal resolution of $M(L)$ by Lemma 3.1 in~\cite{BPS}.
\end{proof}

Proposition~\ref{scarf prop} is illustrated in Figure 2.  The ideal $(cde^2, bde^2, ae^2, a^2bce, a^2bcd)$ has been depolarized as in Theorem~\ref{depolarization prop}; in other words, the polarization $a^2 \mapsto ag$, $e^2 \mapsto ef$ of the given ideal yields $M(L)=(cdef, bdef, aef, abceg, abcdg)$. 
\begin{figure}[h]
\centerline{\epsfig{file=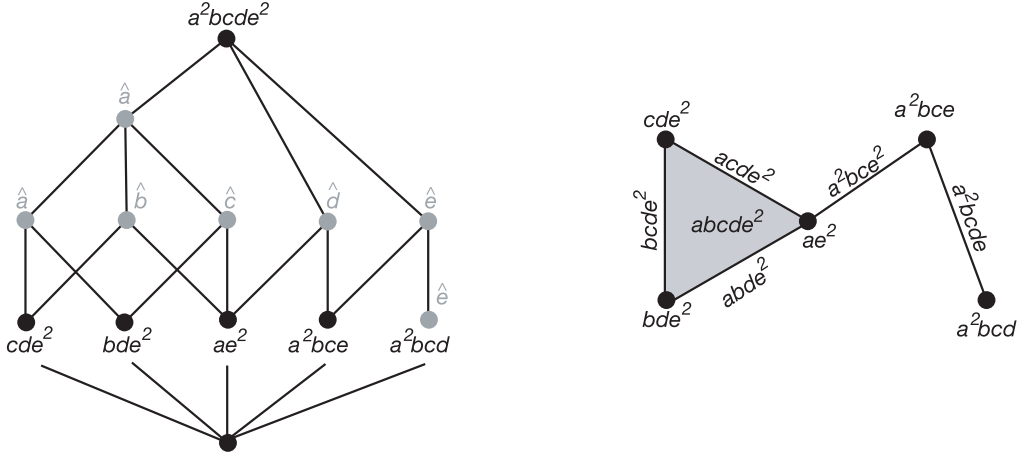 ,height=6cm}}
\caption{The simplicial complex is the Scarf complex of the ideal constructed from its face poset.}
\end{figure}

\begin{definition}\label{minimal definition}\rm
{\em Squarefree minimal ideals} are monomial ideals of the form $M(L)$ for a finite atomic lattice $\hat{L}$.  Say a monomial ideal $I$ is a {\em minimal ideal} iff its polarization is a squarefree minimal ideal.
\end{definition}

The name is justified by Theorem~\ref{universal property theorem}, but note that $M(L)$ is minimal with respect to generators, not with respect to containment.

\begin{theorem}\label{universal property theorem}
Let $I\subset k[{\bf y}] := k[y_1,\ldots,y_n]$ be any squarefree monomial ideal with $LCM(I)\cong \hat{L}$ and let $M=M(L)\subset k[L]$ be the squarefree minimal ideal.
\begin{itemize}
\item[(i)] $n \geq \#\mi{L}$.
\item[(ii)] There exists a $k$-algebra injection $\phi : k[L]\mapsto k[\bf{y}]$ sending variables to variables that satisfies $I\subset \phi(M)\cdot k[\bf{y}]$.  If $n=\#\mi{L}$, then $\phi$ is an isomorphism.
\item[(iii)] Embed $LCM(I)\hookrightarrow 2^n$ and $\hat{L}=LCM(M) \hookrightarrow 2^{\#\mi{L}}$ by identifying a monomial with its support.  Then there exists a join-preserving injection $\rho: 2^{\#\mi{L}} \ra 2^n$ such that $\rho(L)=LCM(I)$.
\end{itemize}
\end{theorem}

We need some definitions before giving the proof.  If $k\in \mi{L}$ and $l\in \hat{L}$ is the unique element covering $k$, then call $(k,l)$ an {\em essential pair} of $\hat{L}$.  Let $I$ be as in the theorem and for $b\in \hat{L}$ let $y(b)\in k[\bf{y}]$ be the monomial corresponding to $b$.  A variable $y_i$ is said to {\em separate} the essential pair $(k,l)$ if $y_i$ divides $y(l)$ but does not divide $y(k)$.  Every essential pair is separated by some variable, but not conversely.

Define for each variable $y_i$ the order filter $L(y_i,1) := \{ b\in \hat{L} :\, y_i \text{ divides } y(b) \}$.  $L(y_i,1)$ is the inverse image of $1$ under the characteristic map $\hat{L} \ra \{0,1\}$ determined by $y_i$.  This map preserves joins so $L(y_i,1)$ is a filter and its complementary order ideal $L(y_i,0) := \hat{L} - L(y_i,1)$ is join-closed.  For $k\in \mi{L}$ define $L(x_k,1)$ and $L(x_k,0)$ analagously.  If $b\in \hat{L}$, then $x_k$ doesn't divide $x(b)$ if and only if $b\leq k$.  Therefore $L(x_k,0)\subset \hat{L}$ equals the principal order ideal $\floor{k}$ generated by $k$.

If $L(y_i,1)=L(y_j,1)$ for some $i\neq j$, then the lcm-lattice cannot distinguish between $y_i$ and $y_j$.  In this sense, one of the variables is unnecessary.

\begin{lemma}\label{filter lemma}
Let $(k,l)$ be an essential pair of $\hat{L}$.  Define $A(k) := \{ i \in [n]: \, L(x_k,1)=L(y_i,1) \}$.  Then $A(k)$ is not empty, and $i\in A(k)$ if and only if $y_i$ separates $(k,l)$. 

Suppose $y_i$ separates no essential pair.  Define $D(i):= \{k\in \mi{L} : \, L(x_k,1)\subsetneq L(y_i,1) \}$.  Then $L(y_i,1)=\cup_{k\in D(i)} L(x_k,1)$.
\end{lemma}
\begin{proof}
Note that $x_k$ separates $(k,l)$, so $i\in A(k)$ implies $y_i$ separates $(k,l)$, too.  Conversely, suppose $y_i$ separates $(k,l)$.  $L(y_i,0)$ is join-closed so it has a maximum element $z$ and $L(y_i,0) = \floor{z}$.  Evidently $\floor{k}\subset \floor{z}$.  If they are not equal, there must exist $k'\in \floor{z} - \floor{k}$.  Then $k'\vee k > k$ implies $k' \vee k \geq l$ implies $k'\vee k \in L(y_i,1)$, a contradiction because $L(y_i,0)$ is join-closed.  Hence $L(y_i,0)=L(x_k,0)$ and this is equivalent to $i\in A(k)$.  $A(k)$ is not empty because $y:L\ra k[{\bf y}]$, $a\mapsto y(a)$, is injective.

For the second statement, assume $y_i$ separates no essential pair.  By the preceding argument, if $L(x_k,1)\subset L(y_i,1)$, then this containment is proper.  Let $z\in \hat{L}$ be the maximum element of $L(y_i,0)$ and let $a\in L(y_i,1)$.  Our task is to find $k\in \mi{L} \cap \ceiling{z}$ such that $k \not \in \ceiling{a}$, for then $a\in L(x_k,1)\subset L(y_i,1)$.  But $a\not\in \floor{z}=L(y_i,0) \iff \ceiling{z}\not\subset \ceiling{a}$.  Equation~\ref{separation equation} implies there exists $k\in (\mi{L}\cap \ceiling{z}) - (\mi{L}\cap \ceiling{a})$ and we are done.
\end{proof}

\phantom{blah}

\begin{proof}{\bf{ of Theorem~\ref{universal property theorem}}}
To prove $(i)$ it suffices to show each $y_i$ separates at most one essential pair.  If $k',k\in \mi{L}$ and $k'\neq k$, then by Lemma~\ref{filter lemma} $A(k)\cap A(k') = \emptyset$ and $(i)$ follows.

For each $k\in\mi{L}$, choose $i(k)\in A(k) \neq \emptyset$.  Define a $k$-algebra morphism $\phi: k[L]\ra k[{\bf y}]$ by $\phi(x_k)=y_{i(k)}$.  Let $a\in \hat{L}$ be an atom and let $x(a)$ and $y(a)$ be the corresponding generators of $M$ and $I$, respectively.  Lemma~\ref{filter lemma} implies $\phi(x(a))$ divides $y(a)$, so $I\subset \phi(M)\cdot k[{\bf y}]$.  This proves $(ii)$.

For each $y_i$ that separates no essential pair, choose a {\em cover} $C(i)\subset D(i)$ such that $L(y_i,1)=\cup_{k\in C(i)} L(x_k,1)$.  The map $\rho$ will depend on the collection of covers ${\cal C} = \{ C( i ) \}$.  Fix ${ \cal C}$.  For $k\in \mi{L}$ define $B(k) = \{i\in [n]:\, k\in C(i)\}$.  Observe that $A(k')\cap B(k) = \emptyset$ for all $k',k\in \mi{L}$, including $k'=k$.

For simplicity regard $2^{\#\mi{L}}$ as the set of subsets of $\mi{L}$.  Define $\rho: 2^{\#\mi{L}} \ra 2^n$ on singletons by $\rho(k):= \rho(\{ k \}) := A(k)\cup B(k) \in 2^n$.  If $k'\neq k$, then $\rho(k')$ and $\rho(k)$ are incomparable because $A(k')\cap (A(k)\cup B(k)) = \emptyset$.  Extend $\rho$ by $\rho(F) = \rho(k_1) \vee \cdots \vee \rho(k_s)$ for each $F=\{k_1,\ldots, k_s\} \in 2^{\#\mi{L}}$.  If $F'\nleq F$, then $A(k)$ is not contained in the support of $\rho(F)$ for any $k\in F' - F$.  This implies $\rho$ is injective.  

To see $\rho(\hat{L}) = LCM(I)$ it suffices to prove $\rho(\supp x(a)) = \supp y(a)$ for all atoms $a\in \hat{L}$.  Recall $\supp y(a) = \{ i\in [n]:\, a\in L(y_i,1)\}$.  Let $k\in \supp x(a)$.  Evidently $A(k)\subset \supp y(a)$.  If $i\in B(k)$, then $a\in L(x_k,1)\subset L(y_i,1)$ so that $i\in \supp y(a)$.  Therefore $\rho(k) \subset \supp y(a)$ and thus $\rho(\supp x(a))\subset \supp y(a)$.  Let $i\in \supp y(a)$.  We can assume $y_i$ separates no variable.  Because $C(i)$ is a cover there exists $k\in C(i)$ such that $a\in L(x_k,1)\subset L(y_i,1)$.  Then $i\in \rho(\supp x(a))$ which implies $\rho(\supp x(a)) = \supp y(a)$.
\end{proof}

\phantom{blah}

For any squarefree ideal $J\subset k[z_1,\ldots,z_m]$, let $SR(J)$ be the simplicial complex on $[m]$ whose Stanley-Reisner ideal equals $J$.

\begin{corollary}
With the identification in Theorem~\ref{universal property theorem} (ii), $SR(M(L))$ is a subcomplex of $SR(I)$.
\end{corollary}

In Section~\ref{inherited properties section} we will show that $SR(M(L))$ inherits many nice properties of $SR(I)$, even though they are generally not homotopy equivalent.  

Note that $k[L] \cong k[{\bf y}]/J$ where $J$ is generated by binomials $y_i-1$.  If $F_{\bullet}$ is a minimal free $k[{\bf y}]$-resolution of $k[{\bf y}]/I$, then $F_{\bullet} \otimes k[{\bf y}]/J$ is a minimal free resolution of $k[L]/M(L)$.  This is an example of the relabeling process in~\cite{GPW}.

We can take this in another direction.  Replacing $\mi{L}$ with $L$ in the construction in Theorem~\ref{main construction theorem} yields a nonminimal ideal $N(L)\subset S(\hat{L}):= k[x_l:\, l\in L]$ whose LCM-lattice equals $L$.  For example, if $\hat{L}$ is the boolean lattice with three atoms $a,b,c$ and coatoms $d,e,f$, then $N(L)=(bcf, ace, abd)$.  This ideal is natural in a different sense:  For any atomic lattice $L$ on $r$ atoms, a (nonminimal) resolution of $N(L)$ can be obtained as a quotient of $N(2^r)$.  To see this, let $a_1,\ldots,a_r$ be the atoms of $L$.  Define the map $\deg:2^r \ra L$ by
$$\deg(F)=\bigvee_{i\in F} a_i.$$
$\deg$ is a join-preserving surjection.  For $l\in L$ the fiber $\deg^{-1}(l)$ is join-closed and therefore has a maximum element.

\begin{theorem}\label{quotient resolution theorem}
Let $\hat{L}$ be a finite atomic lattice on atoms $a_1,\ldots,a_r$.  Let $T$ be the minimal free $S(\hat{L})$-resolution of $N(2^r)$. Let $I\subset S(2^r)$ be the ideal defined by
$$I:=(x_F-1: F\neq \max \deg^{-1}(\deg(F))).$$
Then $S(L)=S(2^r)/I$ and $N(L)=N(2^r)/I\cdot N(2^r)$.  Furthermore $T\otimes S(2^r)/I$ is a (generally nonminimal) free resolution of $N(L)$.
\end{theorem}
\begin{proof}
Let $F$ be a proper subset of $[r]$ and let $x_F\in S(2^r)$ be the corresponding variable.  $x_F-1$ is not a zero divisor on the quotient of a polynomial ring by a monomial ideal, so induction on the number of generators of $I$ shows $T\otimes S(2^r)/I$ is exact.  Fix $l\in L$ and let $A=\max\deg^{-1}(l)\in 2^r$.  Let $F\in \deg^{-1}(l)$ be different than $A$.  It suffices to show $z_A, z_F\in N(2^r)$ have the same image in $S(2^r)/I$.  The quotient of the monomials is $\frac{z_A}{z_F}=\frac{x_{\ceiling{F}}}{x_{\ceiling{A}}}$.  Let $G\in \ceiling{F} - \ceiling{A}$.  Need to see $G\neq \max\deg^{-1}(\deg(G))$.   But $\deg(G\wedge (A-F))=\deg(G)\wedge \deg(A-F)=\deg(G)$ because $\deg(A-F)\leq \deg(A)=\deg(F)\leq \deg{G}$ and $\deg$ preserves order.  $A-F$ is not contained in $G$ (else $G\in \ceiling{A}$) and $A-F$ is nonempty, so $G<G\vee (A-F)$.
\end{proof}

\end{section}

\begin{section}{The Distributive Completion of $\hat{L}$}\label{distributive section}
Recall a finite lattice $J$ is a distributive lattice if for all $a,b,c\in J$
$$a \vee (b \wedge c) = (a \vee b) \wedge (a \vee c).$$
Let $P$ be a finite poset and let $J(P)$ be the set of all order filters in $P$ ordered by reverse inclusion.  $J(P)$ is a distributive lattice and the order dual version of Birkhoff's theorem states that every distributive lattice is of this form: if $\mi{J}$ is the set of meet-irreducible elements of $J-\{\hat{0},\hat{1}\}$, then $J = J(P)$.  See Theorem 3.4.1 in~\cite{Stanley_Enumerative_1}; note that our definition of $J(P)$ is dual to the standard one.

Write $2^{\mi{L}}$ for the set of subsets of $\mi{L}$ ordered by inclusion.  Let $x^c: L \cup \{\hat{1}\} \ra 2^{\mi{L}}$ be the map $x^c(a) := \mi{L}\cap \ceiling{a}_{\hat{L}}$.  The support of $x^c(a)$ is complementary to the support of the monomial $x(a)$ used in Theorem~\ref{main construction theorem}.  $x^c$ is injective by Equation~\ref{separation equation} and sends joins to meets.  

\begin{theorem}\label{distributive theorem}
Let $\phi : 2^{\mi{L}} \ra J(\mi{L})$ be the surjective map $\phi(\{k_1,\ldots,k_s\}) = \ceiling{k_1,\ldots,k_s}$.  $\phi$ is order reversing.  $L\cup \{\hat{1}\}$ is isomorphic as a join semilattice to its image under the composition $\phi\circ x^c$. 
\end{theorem}
\begin{proof}
$\phi$ is order reversing because $J:=J(\mi{L})$ is ordered by reverse inclusion.  Let $\phi'$ be the restriction of $\phi$ to $x^c(L \cup \{\hat{1}\})$.  Then $\phi'$ is injective by the discussion preceding Equation~\ref{separation equation} and it preserves joins because $\ceiling{a\vee b} = \ceiling{a} \cap \ceiling{b}$.  Since $x^c$ is also order reversing, $\phi x^c(L\cup \{\hat{1}\})$ and $L\cup \{\hat{1}\}$ are isomorphic as join semi-lattices.  
\end{proof}

We call $J(\mi{L})$ the distributive completion of $\hat{L}$.  The name is justified because any distributive lattice $J'$ containing $\mi{L}$ as a meet irreducible elements must contain $J(\mi{L})$.  Hence $J'$ must contain $\hat{L}$ as a join subsemilattice.  The theorem also places a restriction on which posets can be isomorphic to $\mi{L}$.

\begin{example}\rm
Using Theorem~\ref{distributive theorem}, one can show the poset $N$ in Figure 3 cannot be isomorphic to $\mi{L}$ for any atomic lattice $L$.  On the other hand, if $\hat{L}=LCM(b^2cd,abd,abc,a^2cd)$, then $\mi{L}$ contains $N$ as a subposet.  See~\cite{Rival} for properties of posets not containing $N$, which are called {\em series-parallel} posets, and their applications in scheduling theory.
\begin{figure}[h]
\centerline{\epsfig{file=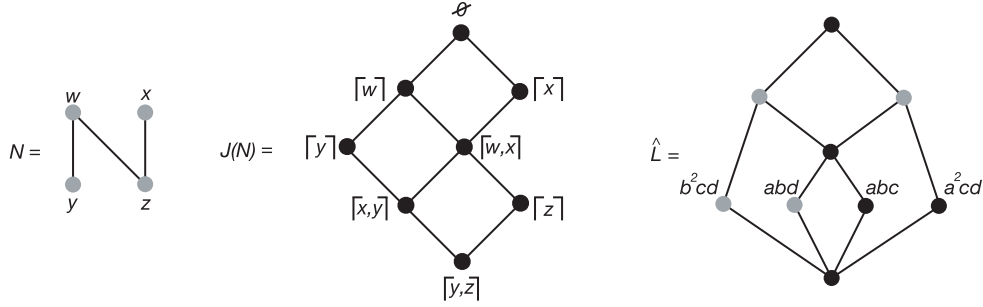 ,height=4cm}}
\caption{$\mi{L}$ can contain $N$ but cannot equal $N$.}
\end{figure}
\end{example}

Let $\phi: 2^{\mi{L}} \ra J(\mi{L})$ be as in Theorem~\ref{distributive theorem}.    Each fiber $\phi^{-1}(F)$ contains a maximum and minimum.  Elements of $\hat{L}=LCM(M(L))$ are the maximums of their individual fibers; associated primes of $M(L)$ are the minimums of their individual fibers.  Furthermore, if we let $P\subset 2^{\mi{L}}$ be the subset of elements $\{k_1,\ldots,k_s\}$ such that the prime ideal $(x_{k_1},\ldots,x_{k_s})$ contains $M(L)$, then $\phi(P)$ is an order ideal of $J(\mi{L})$. It will be useful to have another characterization of the primes containing $M(L)$.  

\begin{prop}\label{primary decomposition prop}
Let $I\subset k[{\bf y}]$ be a squarefree monomial ideal with $\hat{L}\cong LCM(I)$.  The prime $(y_{i_1},\ldots,y_{i_s})\subset k[{\bf y}]$ contains $I$ if and only if
\begin{equation}\label{prime equation}
 \hat{L} - \{\hat{0}\} = \bigcup_{j = 1}^s L(y_{i_j},1) \iff \atom{\hat{L}} \cap \bigcap_{j=1}^s L(y_{i_j},0) \, = \emptyset.
\end{equation}
If this is true, then $(y_{i_1},\ldots, y_{i_s})$ is an associated prime if and only if the union is irredundant if and only if the intersection is irredundant.  
\end{prop}
\begin{proof}
The two conditions in Equation~\ref{prime equation} are easily seen to be equivalent.  $I \subset (y_{i_1},\ldots,y_{i_s})$ if and only if $[n]-F$ is a face of $SR(I)$ where $F=\{i_1,\ldots,i_s\}$.  $[n]-F$ is a face of $SR(I)$ if and only if $\forall$ atoms $a\in \hat{L}\, \exists \, y_{i(a)}\not\in [n]-F$ such that $a\in L(y_{i(a)},1)$ .  This is true if and only if $\hat{L} - \{ \hat{0}\} = \cup_{i\in F} L(y_i,1)$ because the union is an order filter.  Moreover, this union is irredundant if and only if $[n]-F$ is a facet of $SR(I)$.
\end{proof}

\begin{example}\rm
In Figure 4, each element of $J(\mi{L})$ is labeled by the minimal generators of the corresponding order filter.  The image of $L\cup \{ \hat{1}\}$ has been highlighted.  Let $L'$ be the complement of $L$ in $J(\mi{L})-\{\hat{0},\hat{1}\}$.  Note the six associated primes of $M(L)$ appear in $L'\cup \{\hat{0}\}$ under this labeling, as do the two unassociated primes $(a,d,e)$ and $(c,d,e)$.  
$L'$ is not homotopy equivalent to $SR(M(L))$, but it is homotopy equivalent to the subcomplex of $SR(M(L))$ supported on the vertices $c,d,e$, which minimally generate the filter $\hat{0}_J = \{a,b,c,d,e\}$.

Furthermore, $\hat{L}'$ is a coatomic lattice in this case, and therefore it's order dual is a finite atomic lattice.  This is not true in general.  It would be interesting to characterize those atomic lattices $\hat{L}$ for which the order dual of $\hat{L}'$ is also an atomic lattice.
\begin{figure}[h]
\centerline{\epsfig{file=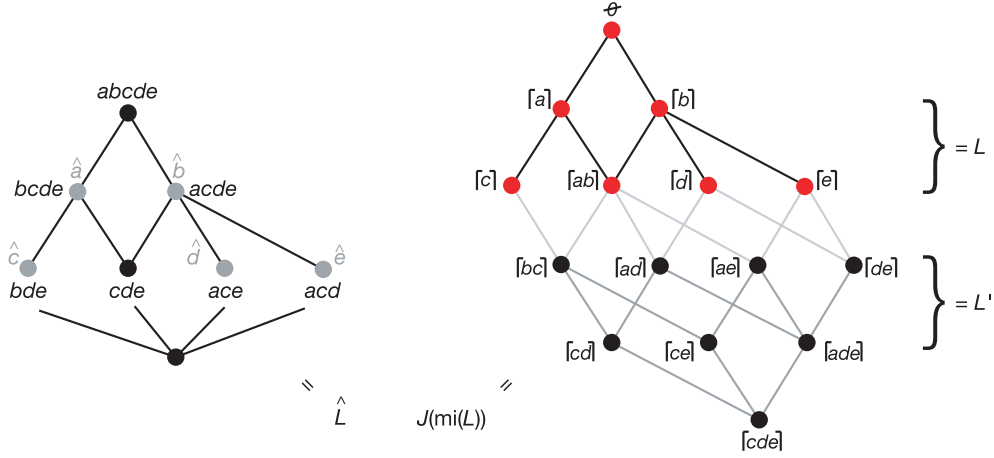, height = 6cm}}
\caption{$(bde,cde,ace,acd)= (b,c)\cap (a,d)\cap (a,e) \cap (d,e) \cap (c,d) \cap (c,e)$.}
\end{figure}
\end{example}

\end{section}

\begin{section}{Nonsquarefree Ideals and Depolarization}\label{depolarization section}
One consequence of Theorem~\ref{betti degree theorem} is that the length of the minimal resolution of $k[{\bf y}]/I$ is $\leq$ to the rank (or height) of $LCM(I)$.   In this section we will characterize which monomial ideals $I\subset k[{\bf y}]$ have polarizations $I_{\pol}$ that are minimal ideals.   As a corollary we will prove that the width of $LCM(I)$ also bounds the projective dimension $\pd k[{\bf y}]/I$ of $k[{\bf y}]/I$ (Theorem~\ref{pd theorem}).

We first recall some facts about the polarization operation.  If $I=(f_1,\ldots,f_r) \subset k[{\bf y}]$ is a monomial ideal, define a simple polarization of $I$ to be the monomial ideal $I'$ obtained as follows.  If $y_i^e$ is the highest power of $y_i$ appearing in any generator $f_1,\ldots, f_r$, then introduce a new variable $y_i'$ and let $I'$ be the ideal obtained by replacing $y_i^e$ by $y_i^{e-1}y_i'$ in each generator $f_j$ divisible by $y_i^e$.  A polarization $I_{\pol}$ of $I$ is a squarefree ideal obtained by iterated simple depolarizations.  If $I_{\pol} \subset S$ where $S$ is a polynomial ring containing $k[{\bf y}]$, then $k[{\bf y}]/I$ is the quotient of $S/{I_{\pol}}$ by an ideal generated by a regular sequence consisting of binomial differences of variables.  Thus many nice properties hold for $I$ if and only if they hold for $I_{\pol}$.  See, for example, \cite{Faridi_polarization} for details on polarization.

\begin{prop}\label{depolarization prop}
Let $\hat{L}$ be a finite atomic lattice and let $M(L)\subset k[L]$ be the minimal ideal.  Let $\gamma_1,\ldots,\gamma_s$ be a partitioning of $\mi{L}$ into $s$ chains.  For each $i$ choose $k(i)\in \gamma_i$.  Define the ideal $J$ to be
$$ J := (x_{l} - x_{k(i)} : \, i\in [s] \text{ and } l \in \gamma_i).$$
Then $k[L]/J$ is a polynomial ring, $\overline{M}:= M(L)\otimes k[L]/J$ is a monomial ideal, and $LCM(\overline{M}) \cong \hat{L}$.
\end{prop}
\begin{proof}
Evidently $\overline{M}_{\pol}$ is combinatorially equivalent to $M(L)$.  It's easy to see $LCM(I_{\pol}) \cong LCM(I)$ for any monomial ideal $I$.
\end{proof}

\begin{corollary}
Let $k',k\in \mi{L}$ and let $\overline{M} = M(L)\otimes k[L]/(x_{k'} - x_k)$.  Then $LCM(\overline{M}) = \hat{L}$ if and only if $k'$ and $k$ are comparable.
\end{corollary}
\begin{proof}
If $k'$ and $k$ are comparable, then $LCM(\overline{M})= \hat{L}$ by the proposition since we can take $\gamma_1 = \{k',k\}$ and let $\gamma_i$ be singletons for $i>1$.  Suppose $k'$ and $k$ are incomparable and let $(k',l')$ and $(k,l)$ be essential pairs.  For all $f\in k[L]$ denote by $\overline{f}$ the image of $f$ in $k[L]/(x_{k'} - x_k)$.  If $l'=l$, then $\overline{x_{k'}} \overline{x_k} = \overline{{x_k}}^2$ divides $\overline{x(l)}$ but divides neither $\overline{x(k')}$ nor $\overline{x(k)}$ and therefore $\lcm (\overline{x(k')},\overline{x(k)}) \neq \overline{x(l)}$.  This implies $LCM(\overline{M})\neq \hat{L}$.  If $l'\neq l$, then $\overline{x_k}$ appears with exponent 1 in $\overline{x(k)}$ and $\overline{x(l)}$ and polarizing cannot distinguish between $\overline{x(k)}$ and $\overline{x(l)}$; therefore $LCM(\overline{M})\neq \hat{L}$
\end{proof}

\begin{corollary}
If $\hat{L}$ is coatomic, then there is no monomial ideal on fewer than $\#\mi{L}$ variables that has lcm-lattice equal to $\hat{L}$.
\end{corollary}
\begin{proof}
If such an ideal $I$ existed, then a minimal one exists  so we can assume $I_{\pol} \cong M(L)$.  $\hat{L}$ is coatomic iff $\mi{L}$ is an antichain, so there would be no way to depolarize $M(L)$ to obtain $I$.
\end{proof}

\phantom{blah}

Recall the width of a poset $P$ is one less than the cardinality of the longest antichain in $P$.

\begin{theorem}\label{pd theorem}
Let $I\subset k[{\bf y}]$ be a monomial ideal with $\hat{L}\cong LCM(I)$.  Then 
$$\pd k[{\bf y}]/I \leq \min(\width \mi{L},\height \hat{L}).$$
\end{theorem}
\begin{proof}
$\pd k[{\bf y}]\leq \height \hat{L}$ by Theorem~\ref{betti degree theorem}. Dilworth's theorem asserts that any poset $P$ can be partitioned into $\width P$ chains.  Proposition~\ref{depolarization prop} implies there exists a monomial ideal $I'\subset S$ in $\width \mi{L}$ variables whose lcm-lattice equals $\hat{L}$.  Theorem~\ref{betti degree theorem} implies $I$ and $I'$ have resolutions of the same length and Hilbert's Syzygy theorem implies $\pd S/I'\leq \width \mi{L}$.  
\end{proof}

Suppose now that $I\subset k[y_1,\ldots, y_n]$ is a squarefree ideal and let $I^a$ be the Stanley-Reisner ideal of the simplicial complex on $[n]$ that is Alexander dual to $SR(I)$.  Terai was the first to prove the connection between the projective dimension of $I$ and the Castenuovo-Mumford regularity of $I^a$.  See also~\cite{BCP}.

\begin{theorem} (Terai~\cite{Terai})
Let $I\subset k[y_1,\ldots, y_n]$ be a squarefree ideal and let $I^a$ be the Alexander dual ideal.  Then
$$\pd k[{\bf y}]/I = \reg I^a.$$
\end{theorem}

The next corollary is a natural consequence.  Since we know there exists a simplicial complex on at most $\#\mi{L}$ vertices that has intersection lattice equal to $\hat{L}$, we also get an easy upper bound on the nonvanishing homology of $SR(I)$.

\begin{corollary}
Let $I$ and $I^a$ be as above.
\begin{itemize}
\item[(i)] $\reg I^a \leq \min(\width \mi{L},\height \hat{L})$.
\item[(ii)] $\tilde{H}_j(SR(I);k) = 0$ if $j< n -\# \mi{L} -1$, where $\tilde{H}_j$ denotes the $j$th reduced homology of $SR(I)$.
\end{itemize}
\end{corollary}
\begin{proof}
$(i)$ follows immediately from Terai's theorem and Theorem~\ref{pd theorem}.  For $(ii)$ we can assume each variable $y_i$ divides some generator of $I$; otherwise SR(I) would be acyclic because its dual would be a cone.  $\hat{L}$ is isomorphic to the intersection lattice of $SR(I)^a$ by Proposition 2.3 in~\cite{GPW} (see also~\cite{Yuzvinsky}).  Therefore the order complex of $L$ is homotopy equivalent to $SR(I)^a$.  The last two sentences are true if $I$ is replaced by $M(L)$, so $(ii)$ follows from Alexander duality.
\end{proof}

\end{section}

\begin{section}{Inherited Properties}\label{inherited properties section}
Henceforth $I\subset k[{\bf y}]$ will be a monomial ideal such that $\hat{L}\cong LCM(I)$.  In this section we prove that the minimal ideal $M(L)$ inherits many nice properties of $I$.  For example, if $I$ has a linear resolution, then so does $M(L)$.  As a corollary, we prove that ideals with linear resolutions are minimal up to a common factor dividing all the generators.  Theorem~\ref{linear res characterization theorem} characterizes those atomic lattices $\hat{L}$ for which $M(L)$ has a linear resolution.  A very similar (by Eagon-Reiner) formulation was proved by Yuzvinsky (\cite{Yuzvinsky}, Theorem 6.4) in the broader setting of rings of sections of sheaves over $L$.

Recall that $I = (m_1,\ldots,m_r)$ is said to have linear quotients if there is an ordering, say the given one, such that $(m_1,\ldots, m_{i-1}): m_i$ is generated by variables for each $i=2,\ldots,r$.  If $I$ is squarefree, then $I$ has linear quotients if and only if $SR(I)^a$
is shellable (\cite{HHZ}, Theorem 1.4).  If $I$ is squarefree, say $I$ is matroidal if $SR(I)$ is a matroid.  $I$ is matroidal if and only if for all $m_i,m_j$ and all $y_l$ dividing $\gcd(m_i,m_j)$, there exists $m_k$ that divides $\lcm(m_i,m_j)/y_l$.

\begin{prop}
Let $I\subset k[{\bf y}]$ be a monomial ideal with $\hat{L}\cong LCM(I)$.
\begin{itemize}
\item[(i)] $\codim I \leq \codim M(L)$.  If $I$ is pure, then so is $M(L)$.
\item[(ii)] $I$ is Cohen-Macaulay  if and only if $\codim I = \codim M(L)$ and $M(L)$ is Cohen-Macaulay.
\item[(iii)] Suppose $I$ is squarefree.  If $I$ has linear quotients then $M(L)$ has linear quotients.  If $I$ is matroidal, then so is $M(L)$.
\item[(iv)] $\reg M(L) \leq \reg I$.
\end{itemize}
\end{prop}
\begin{proof}
Polarizing, if necessary, we can assume $I$ is squarefree.  Then each associated prime of $M(L)$ determines an associated prime of $I$ of the same codimension, by Proposition~\ref{primary decomposition prop}.  This proves $(i)$.  $(ii)$ follows from $(i)$ and the fact that $I$ and $M(L)$ have minimal resolutions of the same length (\cite{GPW}, Theorem 3.3). $(iii)$ is immediate from the preceding definitions.

To prove $(iv)$, embed $k[L]$ in $k[{\bf y}]$ as in Theorem~\ref{universal property theorem}.  The generators of $I$ are obtained from the generators of $M(L)$ through multiplication by new variables.  By induction it suffices to prove the assertion when there is a single new variable $y$.  Set $r=\reg M(L)$ and let $b\in L$ be an $i$th betti degree with $\deg x(b) - i  = r$.  Then $\deg y(b) = \deg x(b)$ or $\deg y(b) = \deg x(b) + 1$, so $\reg M(L) \geq \reg I$. 
\end{proof}

For the rest of the section we will focus on linear resolutions.  Passing to the minimal ideal usually does not commute with Alexander duality, so the next theorem is independent of the previous proposition.

\begin{theorem}\label{linear res theorem}
Let $I\subset k[{\bf y}]$ be a monomial ideal with $\hat{L} \cong LCM(I)$.  If $I$ has a linear resolution, then so does $M(L)$.
\end{theorem}

We postpone the proof.  Lemma~\ref{codim 1 lemma} will be essential for what follows.  Hartshorne proved that every Cohen-Macaulay variety is connected in codimension 1 (see~\cite{Eisenbud_CA}, Theorem 18.12).  Lemma~\ref{codim 1 lemma} asserts that more is true if we limit ourselves to simplicial complexes, which are essentially unions of coordinate planes.  

\begin{lemma}\label{codim 1 lemma}
Let $\Delta$ be a simplicial complex on $[n]$.  Let $\Delta_1\subset \Delta$ be the subcomplex generated by faces $F$ which are the intersections of the facets of $\Delta$ containing $F$.  Assume $\Delta$ is not a simplex.  If $\Delta$ is Cohen-Macaulay, then $\Delta_1$ is Cohen-Macaulay and of codimension 1 in $\Delta$.
\end{lemma}
\begin{proof}
$\Delta = \Delta_1$ if and only if $\Delta$ is a simplex.  Assume $\Delta$ is CM but not a simplex.  If $\dim \Delta = 0$ or 1, then the assertion is easily verified.  By induction, we can assume $\Delta$ is CM implies $\Delta_1$ is CM of codimension 1 if $\dim \Delta < d$.  Let $\dim \Delta = d > 1$.  The Eagon-Reiner theorem along with Theorem~\ref{betti degree theorem} implies $\Delta_1$ is pure of codimension 1.  To see this, let $I\subset k[{\bf y}]$ be the Stanley-Reisner ideal of the Alexander dual $\Delta^a$.  $b\in L:=LCM(I)$ covers an atom if and only if $y(b)$ is the betti degree of a first syzygy because $I$ has a linear resolution.  Call an element that covers an atom a {\em super atom}.  The complement of $y(a)$, $a$ an atom,  is a facet of $\Delta$.  The complement of $y(b)$, $b$ a super atom, is a facet of $\Delta_1$.  Therefore $\Delta_1$ is pure and of codimension 1.

Let $F\in \Delta_1$.  Our task is to show $\link_{\Delta_1} F$ has at most top dimensional homology.  If $F\neq \emptyset$, then $\dim \link_{\Delta} F < d$.  Since $\link_{\Delta} F$ is CM, the induction hypothesis implies $(\link_{\Delta} F)_1$ is CM of codimension 1 in $\link_{\Delta} F$; in particular, it has at most top dimensional homology.  But $(\link_{\Delta} F)_1 = \link_{\Delta_1} F$, so we are done if $F\neq \emptyset$.

It remains to prove $\link_{\Delta_1} \emptyset = \Delta_1$ has no homology in dimension $< d-1$.  $\Delta$ is obtained from $\Delta_1$ by attaching $d$-simplices $G_1,\ldots,G_r$ along codimension 1 faces.  Say $G_i$ is {\em fully attached} if $\partial G_i \subset \Delta_1$.  If $G_i$ is fully attached, then attaching $G_i$ adds only a $d$-cell and cannot affect homology in dimension $< d-1$.  If $G_i$ is not fully attached, then $G_i\cap \Delta_1$ is homeomorphic to a $(d-1)$-disk.  Simplicially collapse $G_i$ along one of its free faces; the remaining complex deformation retracts to $G_i\cap \Delta_1$.  Therefore if $G_i$ is not fully attached, attaching it to $\Delta_1$ cannot affect homology in any dimension.  We've shown that $H_i(\Delta_1;k)=H_i(\Delta;k) = 0$ if $i < d-1 = \dim \Delta_1$.  This concludes the proof.
\end{proof}

\begin{corollary}\label{first betti ideal corollary}
Let $I_1$ be the ideal generated by the betti degrees of the first syzygies of $I$.  If $I$ has a linear resolution, then so does $I_1$.
\end{corollary}
\begin{proof}
We can assume $I$ is squarefree.  Then $\Delta_1 = SR(I_1)^a$. The assertion now follows from Eagon-Reiner and Lemma~\ref{codim 1 lemma}.
\end{proof}

We state the next lemma for easy reference; it is an immediate consequence of Theorem~\ref{betti degree theorem}.

\begin{lemma}\label{localization lemma}
Let $I = (m_1,\ldots,m_r)\subset k[{\bf y}]$ be a squarefree monomial ideal such that $\hat{L} \cong LCM(I)$.  If $b\in \hat{L}$ and $I_{\leq b}:=(m_i:\, m_i \leq y(b))$, then $[\hat{0},b]_{\hat{L}} = LCM(I_{\leq b})$.  Moreover, if $I$ has a linear resolution then so does $I_{\leq b}$.
\end{lemma}

The betti degrees generally do not form an order ideal in $\hat{L}$, not even when $I$ has a linear resolution.  Also, the projective dimension of $I$ can be much smaller than the length of $\hat{L}$ so it is somewhat surprising that we can prove statements about all of $\hat{L}$ when the resolution of $I$ is nicely behaved.

\begin{prop}\label{basic facts prop}
Let $I\subset k[{\bf y}]$ be a monomial ideal with $\hat{L}\cong LCM(I)$.  Suppose $I$ has a linear resolution.
\begin{itemize}
\item[(i)] $\hat{L}$ is a graded lattice.
\item[(ii)] Let $a,b\in \hat{L}$.  If $b$ covers $a\neq \hat{0}$, then $\deg y(b) = \deg y(a) + 1$.
\item[(iii)] If $b\in \hat{L}$ is an $i$th betti degree of $k[{\bf y}]/I$, then $\rank b = i$.
\item[(iv)] Let $b\in \hat{L}$.  Suppose $I$ is squarefree and every $y_i$ divides some generator of $I$.  Then $\length (\hat{0},b)_{\hat{L}}$ $= \dim \link_{SR(I)^a} F$ where $F\in SR(I)^a$ is the face complementary to $y(b)$.  In particular, $\length L$ $= \dim SR(I)^a$.
\end{itemize}
\end{prop}
\begin{proof}
$(i)$:  With the notation of Corollary~\ref{first betti ideal corollary}, $I_1$ also has a linear resolution.  By induction on $\dim SR(I)^a$ we can assume $LCM(I_1)$ is graded, since $\dim SR(I_1)^a = \dim SR(I)^a -1$.  Any maximal chain of $\hat{L}$ intersects $LCM(I_1)=LCM(I)-\atom{L}$ in a maximal chain.  It follows that $L$ is graded and $\length L = \length LCM(I_1) + 1$.

$(ii)$:  By Lemma~\ref{codim 1 lemma} and induction it suffices to prove the statement when $a$ is an atom.  In this case, $b$ is a super atom so that $y(b)$ is a second betti degree of $k[{\bf y}]/I$ and therefore $\deg y(b) = \deg y(a) + 1$ because $I$ has a linear resolution.

$(iii)$:  This follows from $(ii)$ and Theorem~\ref{betti degree theorem}.

$(iv)$:  By Lemma~\ref{localization lemma}, it suffices to prove the statement when $F=\emptyset$.  By induction, we can assume $\length \proper{LCM(I_1)}$ $= \dim SR(I_1)^a$ where $\proper{LCM(I_1)} = LCM(I_1) - \{\hat{0},\hat{1} \}$.  From the proof of $(i)$, $\length L = \length \proper{LCM(I_1)} + 1 = \dim SR(I_1)^a + 1 = \dim (SR(I)^a)_1 +1 = \dim SR(I)^a$.
\end{proof} 

Say a monomial ideal is uniformly generated if its generators all have the same degree.  We need one more result before proving Theorem~\ref{linear res theorem}

\begin{prop}\label{inherited properties prop}
Let $I\subset k[{\bf y}]$ be a monomial ideal with $L\cong LCM(I)$.  If $I$ has a linear resolution, then $M(L)$ is uniformly generated.  Moreover, if $a,b\in \hat{L}$ and $b$ covers $a\neq \hat{0}$, then $\deg x(b) = \deg x(a) +1$.
\end{prop}
\begin{proof}
We will prove the second statement first.  Let $m \in k[L]$ be the product of all $x_k$ such that $k\in \atom{\hat{L}}\cap \mi{L}$.  Set $\hat{L_1} = LCM(I_1)$.  Then $M(L)_1 = m\cdot M(L_1)$ so using Lemma~\ref{codim 1 lemma} we can repeatedly replace $I$ with $I_1$, if necessary, and assume $a$ is an atom.  Evidently $\deg x(b) \geq \deg x(a) +1$.  Because $I$ has a linear resolution, $\deg y(b) = \deg y(a) +1$.  Identify $k[L]$ with as a subring of $k[{\bf y}]$ as in Theorem~\ref{universal property theorem}.  For every atom $a$, $y(a) = x(a) f(a)$ for some monomial $f(a)\in k[{\bf y}]$.  This implies $1 = \deg y(b) - \deg y(a) \geq \deg x(b) - \deg x(a) \geq 1$ and the second claim follows.

In the presence of the second assertion, $M(L)$ is uniformly generated if and only if $L$ is graded.   The first assertion now follows from Proposition~\ref{basic facts prop}.
\end{proof}

\phantom{blah}

\begin{proof}{\bf{of Theorem~\ref{linear res theorem}}}
After polarizing, we can assume $I$ is squarefree.  Modding out by extra variables, we can assume every $y_i$ divides some generator of $I$.  By Eagon-Reiner, it suffices to show $SR(M(L))^a$ is CM.  Let $F\in SR(M(L))^a$.  We can assume $F$ is a face in the intersection lattice, for if not, then $\link F$ is a cone.  We need to prove $\link F$ has at most top dimensional homology.  By Lemma~\ref{localization lemma} we can assume $F=\emptyset$ so that $\link F = SR(M(L))^a$.  Proposition~\ref{basic facts prop} implies $\length L = \dim SR(I)^a$.  Because $SR(I)^a$ and $SR(M(L))^a$ are both homotopy equivalent to $L$ and because $SR(I)^a$ is CM, it suffices to prove $\dim SR(M(L))^a = \dim SR(I)^a$.

Let $\hat{L}_1 := LCM(I_1) = \hat{L} - \atom{\hat{L}}$.  We can assume $\dim SR(M(L_1))^a = \dim SR(I_1)^a$ by induction on the dimension.  If $\atom{\hat{L}}\cap \mi{L} = \emptyset$, then $M(L_1) = M(L)_1$ so that $M(L)_1$ is minimal.  Then $\dim SR(M(L))^a = \dim SR(M(L)_1)^a + 1 = \dim SR(M(L_1))^a + 1 = \dim SR(I_1)^a + 1 = \dim SR(I)^a$.

Now assume $\atom{\hat{L}}\cap \mi{L} = \{ a_1,\ldots, a_s\}$.  Then $x_{a_1}\cdots x_{a_s} \cdot M(L_1) = M(L)_1$.  Proposition~\ref{inherited properties prop} implies $M(L)$ and $M(L)_1$ are uniformly generated, say by generators of degree $d$ and $d+1$, respectively.  Hence the generators of $M(L_1)$ all have degree $d+1-s$.  Then $\dim SR(M(L_1))^a = (\#\mi{L} - s) - (d+1-s)-1 = \#\mi{L} - (d+1) -1 = \dim SR(M(L)_1)^a$.  As above, this implies $\dim SR(M(L))^a = \dim SR(I)^a$ and the proof is complete.
\end{proof}

We close this section with a characterization of which finite atomic lattices support a linear resolution.  Yuzvinsky proved a similar statement in the setting of rings of sections of sheaves on poets (Theorem 6.4, \cite{Yuzvinsky}).   As a corollary, we prove that any monomial ideal $I$ that has a linear resolution is essentially minimal.

\begin{theorem}\label{linear res characterization theorem}
Let $\hat{L}$ be a finite atomic lattice.  The minimal ideal $M(L)$ has a linear resolution if and only if 
\begin{itemize}
\item[(i)] $L$ is Cohen-Macaulay as a poset and 
\item[(ii)] if $a,b\in \hat{L}$ and $b$ covers $a\neq \hat{0}$, then $\deg x(b)=\deg x(a) + 1$.
\end{itemize}
\end{theorem}
\begin{proof}
Suppose $(i)$ and $(ii)$ hold.  In the presence of $(ii)$,  $M(L)$ is uniformly generated if and only if $L$ is graded.  Let $d$ be the degree of the generators of $I$.  Then $(ii)$ implies $\length(\hat{0},b)_L = \deg x(b) - d - 1$.  Because $L$ is CM, Theorem~\ref{betti degree theorem} implies the $\beta_i(k[L]/M(L),m)$ can only be nonzero if $i-2=\deg x(b) - d -1$, or equivalently $i=\deg x(b) -d +1$.  This implies $k[L]/M(L)$, and thus $M(L)$, has a linear resolution.

Conversely, suppose $M=M(L)$ has a linear resolution.  Statement $(ii)$ follows from Proposition~\ref{basic facts prop}.  Let $\Delta = SR(M)^a$.  Now consider $\hat{L}$ as the intersection lattice of $\Delta$.  To keep the argument uniform we will work with the polar face poset $P=P(\Delta)^*$, which naturally contains $\hat{L}$ as a subposet.  Let $F,G\in \hat{L}$.  We need to show the order complex of $(F,G)_{\hat{L}}$ has at most top dimensional homology.  By Lemma~\ref{localization lemma} we can assume $G=\hat{1}$.  If $F=\hat{0}$, then $(F,\hat{1})_{\hat{L}}=L$ is homotopy equivalent and of equal dimension to $\Delta$.  Since $\Delta$ is CM, $\proper{L}$ has at most top dimensional homology.  If $F\neq \hat{0}$, then repeatedly replacing $M$ by $M_1$ if necessary, we can assume $F$ is an atom of $\hat{L}$, i.e. a facet of $\Delta$.  

Proposition~\ref{basic facts prop} implies $\length (F,\hat{1})_{\hat{L}} = \dim \Delta - 1$.  Let $H_1,\ldots,H_s\in L$ be the super atoms covering $F$.  Then
$$(F,\hat{1})_{\hat{L}} = \ceiling{H_1,\ldots,H_s}_{\hat{L}_1} - \{\hat{1}\}$$
where $\hat{L}_1 = \hat{L}-\atom{\hat{L}}$ is the LCM-lattice of $M_1$.  Let $H$ be the intersection of the $H_i$.  The proper part of any poset is homotopy equivalent to the proper part of the join closure of its atoms, so $\ceiling{H_1,\ldots,H_s}_{\hat{L}_1}-\{\hat{1}\}$ is homotopy equivalent to $(\hat{0},H)_{\hat{L}_1}$.

Eagon-Reiner imply $\Delta$ is CM.  Let $X\subset \Delta$ be the codimension 1 subcomplex of $\Delta$ along which $F$ is attached.  The facets of $X$ are $H_1,\ldots,H_s$.  Evidently $\ceiling{H_1,\ldots,H_s}_P-\{\hat{1}\}$ equals the proper part of $P(X)^*$.  If $F$ is fully attached, then $X=\partial F$ is a sphere of dimension $\dim \Delta -1 = \length (F,\hat{1})_L$.  If $F$ is not fully attached, then $X$ is a codimension 1 disk, and hence is contractible.  In any case, $(F,\hat{1})_P$ has no homology in dimension $< \length (F,\hat{1})_L$.  Since $\ceiling{H_1,\ldots,H_s} - \{ \hat{1}\}$ is also homotopy equivalent to $(\hat{0},H)_{L_1}$, the proof is complete.
\end{proof}

\phantom{blah}

We can now show that monomial ideals that have linear resolutions are essentially minimal ideals.

\begin{corollary}\label{linear res structure corollary}
Let $I=(m_1,\ldots,m_r) \subset k[{\bf y}]$ be a squarefree monomial ideal.  Assume that $\hat{L}$ $=LCM(I)$ and that $I$ has a linear resolution.  If $\gcd\{m_1,\ldots,m_r\} = 1$, then $I\cong M(L)$.
\end{corollary}
\begin{proof}
Embed $k[L]$ in $k[{\bf y}]$ as in Theorem~\ref{universal property theorem}.  Suppose $L(y_i,1)\neq L -\{\hat{0}\}$ for some $y_i$ that doesn't separate an essential pair.  Let $a\in L$ be a maximal element not in $L(y_i,1)$ and suppose $b\in L(y_i,1)$ covers $a$.  By Lemma~\ref{inherited properties prop}, $\deg x(b) = \deg x(a) +1$.  Multiplying the relevant generators by $y_i$ increases the degree of $b$ but not the degree of $a$.  Hence $\deg y(b) -\deg y(a) \geq 2$, contradicting Propostion~\ref{basic facts prop}.  This implies $L(y_i,1)=L-\{\hat{0}\}$ if $y_i$ doesn't separate an essential pair.

Since $\gcd\{m_1,\ldots,m_r\}=1$, every $y_i$ must separate an essential pair.  An argument similar to the one above shows no two $y_i$ can separate the same essential pair, i.e. $\# A(l)=1$ for all $l\in \mi{L}$.  Modding out by extra variables we can assume every $y_i$ divides some generator of $I$.  Hence $n=\# \mi{L}$.  This implies that $\rho: 2^{\#\mi{L}} \ra 2^n$ is an isomorphism of lattices, and the claim follows.
\end{proof}

\phantom{blah}

\noindent{\bf Acknowledgements.} The author thanks Dave Bayer for many helpful discussions and suggestions.  

\end{section}

\bibliographystyle{plain} 
\bibliography{thesis}

\end{document}